\documentclass[11pt]{article}

\usepackage{amsmath,amsfonts,amssymb,amsthm,mathtools,bbm, geometry}

\DeclareMathAlphabet{\mathboondoxfrak}{U}{BOONDOX-frak}{m}{n}

\usepackage{amssymb}

\newtheorem{thm}{Theorem}[section]

\newtheorem{dfn}[thm]{Definition}

\newtheorem{lem}[thm]{Lemma}

\newtheorem{prop}[thm]{Proposition}

\newtheorem{remark}[thm]{Remark}

\newtheorem{cor}[thm]{Corollary}

\newtheorem{ex}[thm]{Example}

\newtheorem{question}[thm]{Question}

\def\sq{{\scriptscriptstyle \square}}

\def\bq{\begin{question}}
\def\bt{\begin{thm}}
\def\bp{\begin{prop}}
\def\blem{\begin{lem}}
\def\bd{\begin{dfn}}
\def\br{\begin{remark}}
\def\bc{\begin{cor}}
\def\bex{\begin{ex}}
\def\beqs{\begin{eqnarray*}}
\def\beq{\begin{eqnarray}}
\def\bi{\begin{itemize}}

\def\eq{\end{question}}
\def\et{\end{thm}}
\def\ep{\end{prop}}
\def\elem{\end{lem}}
\def\ed{\end{dfn}}
\def\er{\end{remark}}
\def\ec{\end{cor}}
\def\eex{\end{ex}}
\def\eeqs{\end{eqnarray*}}
\def\eeq{\end{eqnarray}}
\def\ei{\end{itemize}}

\def\N{{\mathbb{N}}}

\def\c{\cdot}

\def\ov{\overline}

\def\B{{\cal{B}}} 
\def\fB{{\mathfrak B}}

\def\w*{w^*-w^*}

\def\ra{\rightarrow}

\def\bs{\backslash}

\def\Q{\mathbb{Q}}

\def\fB{{\mathboondoxfrak B}}

\def\nuB{\nu_{\mathcal B}}
\def\nub{\nu_{\mathcal B}}

\def\nuap{\nu_{AP}}

\def\nuwap{\nu_{WAP}}

\def\nulmc{\nu_{LMC}}

\def\Bluc{\B^{LUC}}
\def\Bap{\B^{AP}}
\def\Bwap{\B^{WAP}}
\def\Blmc{\B^{LMC}}

\def\U{{\mathcal U}}
\def\V{{\mathcal V}}

\makeatletter
\def\moverlay{\mathpalette\mov@rlay}
\def\mov@rlay#1#2{\leavevmode\vtop{%
   \baselineskip\z@skip \lineskiplimit-\maxdimen
   \ialign{\hfil$\m@th#1##$\hfil\cr#2\crcr}}}
\newcommand{\charfusion}[3][\mathord]{
    #1{\ifx#1\mathop\vphantom{#2}\fi
        \mathpalette\mov@rlay{#2\cr#3}
      }
    \ifx#1\mathop\expandafter\displaylimits\fi}
\makeatother

\usepackage[all]{xy}

\def\R{\mathbb{R}}

\begin{document}

\title{Totally Disconnected Semigroup Compactifications: Non-Introversion of the Full Boolean Algebra of Clopen Sets}
\author{Joshua Basman Monterrubio\footnote{This undergraduate student author was supported by R. Stokke's NSERC grant in the summer of 2022.},  Thomas Czyzowicz\footnote{This undergraduate student author was  supported by a 2024 NSERC USRA.}, Ross Stokke\footnote{
This author was partially supported   by an NSERC grant.}, and Emily Thevenot\footnote{This undergraduate student  author was supported by R. Stokke's NSERC grant in the summer of 2025.}    }
\date{}
\maketitle

\begin{abstract}{\small  In terms of  the existence of a single clopen set and two related nets, we characterize when the full Boolean algebra, $\fB(G)$,  of clopen subsets of a topological group $G$  is left introverted. We employ this characterization to show that when $G$ is   a first countable,  $\sigma$-compact, totally disconnected locally compact group, $\fB(G)$ is left introverted if and only if $G$ is compact or discrete,  thus providing a strong  positive answer to a question posed in \cite{Ste-Sto}. Examples of clopen sets and nets witnessing our non-introversion theorem are presented. Some hereditary properties of left introversion of $\fB(G)$ are proved and then employed to extend our main result to  other classes of topological groups.   

\smallskip

\noindent{\em MSC codes:}  Primary 22A10, 22A20, 06E15, 54D35 \\
{\em Key words and phrases:} semigroup compactification, totally disconnected locally compact group, Boolean algebra, Stone space, introversion,  Arens product}
\end{abstract} 

\section{Introduction } 

 Totally disconnected locally compact (tdlc) groups are an important class of topological groups that have been studied by many authors in a wide variety of  contexts.  A  recent survey article on tdlc groups is \cite{Cap-Wil}; also see \cite{Bur, Cap-Mon, Pal, Wil}.  Semigroup compactification theory has meanwhile served as  a fundamental tool   in the study of topological groups and semigroups, e.g., see \cite{Ber-Jun-Mil, Rup}. Motivated by the desire to remain within the class of totally disconnected spaces when studying the semigroup compactifications
of a totally disconnected group, the authors of  \cite{Ste-Sto}  introduced definitions of left/right introversion of a Boolean algebra, $\B$,  of clopen subsets of a topological group $G$ and left/right Arens products on the associated  Stone space,  $\nu_\B G$, with respect to which $\nu_\B G$   is a totally disconnected (right/left topological) semigroup compactification of $G$. Moreover, every totally disconnected semigroup compactification of $G$ is equivalent to $\nu_\B G$, with left Arens product, for a unique left introverted Boolean algebra of clopen subsets of $G$ \cite[Theorem 3.3]{Ste-Sto}. The universal totally disconnected right topological semigroup, universal  semitopological semigroup, and universal topological group compactifications of $G$, $\nulmc G$, $\nuwap G$ and $\nuap G$, respectively,  were described in \cite{Ste-Sto}  by their associated Boolean algebras of clopen subsets of $G$, $\Blmc$, $\Bwap$ and $\Bap$, respectively.  

 To help motivate the  theory developed in \cite{Ste-Sto}, the referee of that paper suggested that the authors should aim to  find an example of a topological group $G$ for which the universal totally disconnected compactification of $G$, $\nu_\fB G$, is not equivalent to the universal semigroup compactification of $G$, $\nulmc G$ \cite[Theorem 4.1]{Ste-Sto}.  The referee was thus asking for an example of $G$ for which  $\fB$,  the full Boolean algebra of clopen subsets of $G$, is  strictly larger than $\Blmc$; equivalently, the referee was looking for an example of $G$ for which  $\fB$ is not left introverted. Not knowing of an example, the authors posed this as part of Question 3 on p. 324 of \cite{Ste-Sto}.  The purpose of this note is to provide a strong positive answer to this question: among other results, we will show that whenever $G$ is a first countable (equivalently metrizable), $\sigma$-compact tdlc group, $\fB$ is not left introverted when $G$ is nondiscrete and noncompact (and conversely).  It should be noted, however, that if $G$ is any tdlc group (more generally a non-Archimedean topological group), the totally disconnected semigroup compactifications determined by the left introverted Boolean algebras of clopen sets, $\Omega$ (the open coset ring of $G$), $\Bwap$, $\Bluc$ and $\Blmc$, as defined in \cite{Ste-Sto},  are all topologist's compactifications in the sense that their associated compactification homomorphisms are homeomorphisms onto their images \cite[Corollary 4.15]{Ste-Sto}; that is, though they may be properly contained in $\fB$,   each of  the left introverted Boolean algebras $\Omega$, $\Bwap$, $\Bluc$ and $\Blmc$ contain sufficiently many clopen sets to yield topologist's compactifications of tdlc groups.

 Let $X$ be a (always Hausdorff)  topological space. We let $\fB(X)$, or just $\fB$, denote the full Boolean algebra of all clopen subsets of $X$.  If $\B$ is any Boolean algebra of clopen subsets of $X$, $\nub X$ is its Stone space: this is the set of all $\B$-ultrafilters, which becomes a compact, zero-dimensional Hausdorff space with respect to the topology determined by the sets $$\lambda_\B(B)= \{ \U \in \nub X: B \in \U\} \quad (B \in \B)$$ as an open basis. By the Stone representation theorem, $\{ \lambda_\B(B) : B \in \B\}$ is precisely the set of clopen subsets of $\nub X$ and $B \mapsto \lambda_\B(B)$ is a Boolean algebra isomorphism of $\B$ onto $\fB(\nub X)$. Letting $\nub(x)$ denote the fixed $\B$-ultrafilter at $x$ in $X$, the pair $(\nub X, \nub)$ is a zero-dimensional compactification of $X$, and is a topologist's compactification  of $X$ exactly when $X$ is zero-dimensional and $\B$ is an open basis for $X$. Every zero-dimensional compactification of $X$ is equivalent to $(\nub X, \nub)$ for a unique Boolean algebra of clopen sets $\B$ \cite[Proposition 2.3]{Ste-Sto}. 
 For $x \in X$ and $B \in \B$, it is helpful to note that 
$$x \in B, \  B \in \nuB(x) \text{ and } \nuB(x) \in \lambda_\B(B) \text{ are equivalent statements.}$$ 
 The reader is referred to Sections 3.2 and 4.7 of \cite{Por-Woo}, and Section 2 of \cite{Ste-Sto}, for additional definitions and details.  
 
Assume now that    $G$ is a (Hausdorff) topological group. For a subset $B$ of $G$ and $s \in G$,  let $B \c s = s^{-1}B$ and  $ s\c B = Bs^{-1}$. If $\B$ is  a $G$-invariant Boolean subalgebra of $\fB(G)$ and $\V\in \nub G$, we let
$$\V \c B = \{ s \in G: B\c s \in \V \} \quad \text{ and } \quad B \c \V = \{ s \in G: s \c B \in \V\},$$
and say that $\B$  is \it left \rm  (respectively \it right\rm) \it introverted \rm  if 
$\V \c B \in \B$ (respectively $B \c \V \in \B$) for every $\V \in \nub G$ and $B \in \B$. 
When $\B$ is left (respectively right) introverted, $(\nuB G, \nuB)$ becomes a totally disconnected right (left) topological semigroup compactification of $G$ with respect its left (right) Arens product defined by 
$$\U \sq \V = \{ B \in \B: \V \c B \in \U\} \quad \text{(respectively } \U \diamond \V = \{ B \in \B: B \c \U \in \V\}).$$ Moreover, every totally disconnected right  topological semigroup compactification of $G$ is equivalent, as a semigroup compactification, to $(\nuB G, \nuB)$ for some left introverted Boolean subalgebra $\B$ of $\fB(G)$ \cite[Theorem 3.3]{Ste-Sto}.  Additional details, definitions,  references, and motivation for these notions can be  found in \cite{Ste-Sto}.    
   
   \section{Non-introversion theorems and examples}

   \bp \rm  \label{General Non-Introversion Prop}  Suppose that $B$ is a clopen subset of a topological group $G$ for which \bi \item[(N$_B$):]  there are nets $(a_\gamma)_{\gamma \in \Gamma}$ in $B$ and $(c_\kappa)_{\kappa \in K}$ in $G$ such that $\lim c_\kappa = e_G$, and for each $\kappa$  in $K$ there is some $\gamma_\kappa$ in $\Gamma$ such that $c_\kappa a_\gamma \notin B$ for $\gamma \geq \gamma_\kappa$.  \ei 
   If $\B$ is any $G$-invariant Boolean algebra of clopen sets containing $B$, then there exists some $\V$ in $\nuB G$ such that $\V \c B$ is not open. Thus, any Boolean algebra of clopen sets containing $B$ is not left introverted.   
   \ep 
   
   \begin{proof}  Let $\B$ be a $G$-invariant Boolean algebra of clopen sets that contains $B$. Since $\nuB G$ is compact, there is a subnet $(\nuB(a_{\gamma_i}))_{i \in I}$ of $(\nuB(a_\gamma))_\gamma$ and $\V$ in $\nuB G$ such that $\lim_i  \nuB(a_{\gamma_i}) =\V$ in $\nuB G$. Each $a_{\gamma_i}$ is in $B$, equivalently each $\nuB(a_{\gamma_i})$ belongs to the closed set $\lambda_\B(B)$, so $\V \in \lambda_\B(B)$. Hence, $B = B\c e_G \in \V$, which means that  $e_G = \lim_\kappa c_\kappa \in \V \c B$. Supposing that $\V \c B$ is open, some  $c_\kappa$ belongs to $ \V \c B$. Equivalently $B \c c_\kappa \in \V$,  and therefore $\V \in \lambda_\B(B\c c_\kappa)$. As $\lim_i \nuB(a_{\gamma_i}) = \V$, there is some $i_0$ in $I$ such that $\nuB(a_{\gamma_i}) \in \lambda_\B(B \c c_\kappa)$ for $i \geq i_0$. Hence, $a_{\gamma_i} \in B \c c_\kappa = c_\kappa^{-1}B$, and therefore $c_\kappa a_{\gamma_i} \in B$ for $i \geq i_0$, a contradiction. We conclude that $\V \c B$ is not an open subset of $G$. 
    \end{proof}

\bt \rm   \label{Non-Introversion Characterization Thm} Let $G$ be a topological group. Then  
 $\fB$ is \it not \rm left introverted if and only if 
 there is a clopen subset $B$ of $G$ with  property (N$_B$).   \et 

\begin{proof}  Supposing that  $\fB$ is not left introverted, there is a   $\fB$-ultrafilter $\V$ and a clopen subset $A$ of $G$ such that $\V \c A$ is not clopen; since $G \bs (\V \c A) = \V\c (G \bs A)$, we can assume that $\V \c A$ is not closed.  We can also assume that  $e_G  \in \ov{\V \c A} \bs \V\c A$. (To see this, take $x \in  \ov{\V \c A} \bs \V\c A$. Since $g \mapsto g \c x :=x^{-1}g$ is a homeomorphism of $G$, $e_G  \in  (\ov{\V \c A} \bs \V\c A) \c x =  \ov{(\V \c A)\c x} \bs (\V\c A) \c x =  \ov{\V \c (A\c x)} \bs \V\c (A\c x)$, so we can replace $A$ by $A \c x$ if necessary.)  Take a net $(c_\kappa)_{\kappa \in K}$ in $\V \c A$ such that $\lim_\kappa c_\kappa = e_G$ and let $B = G \bs A$. Since $e_G$ is not in $\V \c A$, $A = A\c e_G$ does not belong to the $\fB$-ultrafilter $\V$, so $B \in \V$; equivalently, $\V \in \lambda_\fB (B)$. Let $(a_\gamma)_{\gamma \in \Gamma}$ be a net in $G$ such that $\lim_\gamma \nu_\fB(a_\gamma) = \V$. Assuming without loss of generality that each $\nu_\fB(a_\gamma)$ belongs to the open neighbourhood $\lambda_\fB(B)$ of $\V$,  $(a_\gamma)_{\gamma \in \Gamma}$ is a net  contained in $B$. Fix $\kappa$ in $K$. Since $c_\kappa \in \V \c A$, $A\c c_\kappa \in \V$, and therefore $\lambda_\fB(A \c c_\kappa)$ is an open neighbourhood of $\V$ in $\nu_\fB G$; hence, we can take $\gamma_\kappa$ in $\Gamma$ such that $\nu_\fB(a_\gamma) \in \lambda_\fB(A \c c_\kappa)$ for $\gamma \geq \gamma_\kappa$. Equivalently, $c_\kappa a_\gamma \in A$, and therefore $ c_\kappa a_\gamma \notin B$, for $\gamma \geq \gamma_\kappa$. Hence $B$ has property (N$_B$). Proposition \ref{General Non-Introversion Prop} includes the converse. 
 \end{proof}

\bp  \rm If $G$ is either a discrete group or a compact group, then $\fB$ is  introverted.    \ep 

\begin{proof} Trivially, $\fB$ is introverted when $G$ is discrete. For $s \in G$ and $B \in \fB$, $\nu_\fB(s) \c B = s \c B, B \c \nu_\fB(s) = B\c s \in \fB$ \cite[Lemma 3.2]{Ste-Sto} and every $\fB$-ultrafilter is fixed when $G$ is compact. Hence,  $\fB$ is introverted when $G$ is compact.  \end{proof} 

\bt  \label{MainThm - metrizable case} \rm Let $G$ be a first countable, $\sigma$-compact  tdlc group. Then $\fB$ is left introverted if and only if  $G$ is discrete or compact.   \et

\begin{proof} Suppose that $G$ is nondiscrete and noncompact.   Let $(C_n)$ be a local basis at the identity composed  of compact open subgroups of $G$ (e.g., see \cite[Theorem 7.8]{HR}); since $G$ is not discrete, we can assume that each $C_{n+1}$ is properly contained in $C_n$.  Every compact set is covered by finitely many cosets of the open subgroup $C_1$, and $G$ is $\sigma$-compact and noncompact, so $C_1$ has countable, infinite index in $G$.  Let $\{a_n: n \in \N \}$ be a complete set of representatives of right cosets of $C_1$, without redundancies.  Observe that because $G$ is noncompact, $\lim_n a_n = \infty$, (meaning that for any compact subset $K$ of $G$, $a_n \notin K$ for large enough $n$). 
Let  $$B_n:= C_n a_n \quad \text{ and } \quad  B := \bigcup_{n\in \N} B_n,$$ an open subset of $G$.   To see that $B$ is also closed, take $y$ in $\overline{B}$ and a sequence $(y_k)_k$  in $B$ that converges to $y$. Suppose that for each $N$, $\{y_k: k \in \N\}$ is \it not \rm contained in $\bigcup_{n=1}^N B_n$. We can then recursively build subsequences $(y_{k_l})_l$ of $(y_k)_k$ and $(B_{n_l})$ of $(B_n)_n$ such that $y_{k_l} \in B_{n_l} = C_{n_l}a_{n_l}$ for each $l$ in $\N$.  But then  $y_{k_l} a_{n_l}^{-1} \in C_{n_l}$  for each $l$, so $\lim_l y_{k_l} a_{n_l}^{-1} = e_G$; therefore $\lim_l a_{n_l} = \lim_l (y_{k_l} a_{n_l}^{-1})^{-1} y_{k_l} = e_G \c  y = y$, a contradiction because $\lim_l  a_{n_l} = \infty$. We conclude that the sequence $(y_k)_k$ is contained in the closed set $\bigcup_{n=1}^N B_n$ for some $N$, so $y = \lim y_k$ belongs to $\bigcup_{n=1}^N B_n$, and therefore belongs to $B$. Hence, $B$ is clopen. 

For each $k$, take $c_k \in C_k \bs C_{k+1}$, so $\lim c_k = e_G$. To see that $B$ satisfies (N$_B$), fix $k$ and suppose that $n \geq k+1$.    Since $c_k \notin C_n$,  $c_ka_n \notin C_n a_n = B_n$.  For $m \neq n$, $C_1 a_n$ and $C_1a_m$  are disjoint, so $C_k a_n$ and $C_m a_m = B_m$ are disjoint; thus, $c_k a_n \notin B_m$.  As $c_ka_n \notin  \bigcup_{m \in \N} B_m = B$, $B$ satisfies property (N$_B$). By Theorem \ref{Non-Introversion Characterization Thm},  $\fB$ is not left introverted.   \end{proof} 

\bex \rm1.  By Theorem \ref{MainThm - metrizable case},   $\fB(\mathbb{Q}_p)$ fails to be left introverted, where $\Q_p$ is the tdlc group of $p$-adic numbers under addition.  Via a slightly simpler argument,  we offer a concrete construction of  a clopen set $B$ satisfying (N$_B$) in this case: Observe that if  $x \in B(p^{-n}; 1)$, the (cl)open ball with radius 1 and centre  $p^{-n}$ for some $n  \in \N$, then $|x|_p = p^n$, (because otherwise $|x-p^{-n} |_p = \max\{ |x|_p, |p^{-n}|_p\} \geq p^n >1$  \cite[Proposition 2.3.4]{Gou}). For each $n \in \N$, let $B_n= B(p^{-n};p^{-n})$, a clopen subset of $B(p^{-n};1)$, and let $B = \bigcup_{n \in \N} B_n$. To see that $B$ is closed, hence clopen, let $(y_k)$ be a sequence in $B$ converging to $y$ in $\Q_p$. Then $(|y_k|_p)$ is a convergent sequence in the discrete set $\{ p^n: n \in \N\}$, hence eventually constant. We conclude that $(y_k)_k$ is eventually in the closed set $B_{n_0}$ for some $n_0$, so $y \in B$. Hence,  $B$ is closed. Taking $c_k = p^k$ and $a_n = p^{-n}$, $\lim c_k =0$ and $(a_n)$ is contained in $B$.  For $k \in\N$, take $n \geq k+1$. Since $p^k \notin B(0; p^{-n})$, $p^k +p^{-n} \notin p^{-n}+B(0;p^{-n}) = B_n$; if $m \neq n$, then $p^k + p^{-n} \notin B_m $ because otherwise $p^k + p^{-n}$ belongs to the disjoint sets $B(p^{-n};1)$ and $B(p^{-m};1)$. Hence $c_k + a_n \notin B$. 

\medskip 

\noindent 2. The set of rational numbers $\Q$, under addition and relative Euclidean topology, is an example of a  non-locally compact totally disconnected topological group $G$  such that $\fB(G)$ is not left introverted: For every positive integer $n$, let $d_n$ be a positive irrational number less than $1/2n$ and let $B = \bigcup_{n \geq 1} (n-d_n, n+d_n)$. Then $B$ is an open subset of $\Q$ and it is closed because $$\Q \bs B = (-\infty, 1-d_1) \cup  \bigcup_{n \geq 1} (n+d_n, n+1-d_{n+1})$$ is  open. To see that $B$ satisfies (N$_B$), let $c_k = 1/2k$ and $a_n = n$. Take $k \in \N$ and suppose that $n\geq k$. Observe that 
$$  n +{1\over 2k} \in \left[n+{1\over 2n}, n+ {1\over 2 } \right]  \subseteq \left[n+{1\over 2n}, n+1 - {1\over 2n+1} \right] \subseteq [n+d_n, n+1- d_{n+1} ] \subseteq \Q \bs B,$$ 
so $a_n + c_k \notin B$. 
 \eex

\bp \label{Prop LI Quotients Open Subgroups}  \rm Let $N$ be a closed normal subgroup of a locally compact group $G$, and let $H$ be an open subgroup of $G$. If $\fB(G)$ is left introverted, then  $\fB(G/N)$ and $\fB(H)$ are left introverted. 
\ep 

\begin{proof}  Suppose that $\fB(G/N)$ is not left introverted. By Theorem \ref{Non-Introversion Characterization Thm}, $G/N$ contains a clopen subset $B$ that satisfies property (N$_B$); take nets $(a_\gamma)_{\gamma \in \Gamma}$ and $(c_\kappa)_{\kappa \in K}$ in $G$ such that  $(a_\gamma N)_{\gamma \in \Gamma}$  is contained in $B$, $\lim_\kappa c_\kappa N = e_{G/N}$, and for each $\kappa \in K$ there is some $\gamma_\kappa$ in $\Gamma$ such that $c_\kappa a_\gamma N \notin B$ for $\gamma \geq \gamma_\kappa$.  The quotient homomorphism $q_N: G \ra G/N$ is an open, continuous map, so we can assume that $(c_\kappa)$ is contained in a compact neighbourhood of the identity; passing to a subnet, if necessary, we can further assume that $(c_\kappa)$ converges to some $c$ in $G$.  Since $cN = \lim c_\kappa N= e_{G/N}$, $c \in N$, and therefore $c_\kappa c^{-1} N = c_\kappa N$; if needed,  by replacing $c_\kappa$ with $c_\kappa c^{-1}$ we can thus assume that $\lim c_\kappa = e_G$. Letting $B_0 = q_N^{-1}(B)$, a clopen subset of $G$, using the nets $(a_\gamma)_\gamma$ and $(c_\kappa)_\kappa$, one sees that $B_0$ satisfies property (N$_{B_0}$). By Theorem \ref{Non-Introversion Characterization Thm}, $\fB(G)$ is not left introverted. If $\fB(H)$ is not left introverted, it is clear from Theorem \ref{Non-Introversion Characterization Thm} that $\fB(G)$ is not left introverted.  
\end{proof} 

The following lemma will be known, but we could not find a reference. 

\blem \label{Lemma Sigma-compact open subgroup}  \rm  If $G$ is a noncompact locally compact group, then it contains a noncompact, $\sigma$-compact open subgroup. 
\elem

\begin{proof}  If $G$ has a compact symmetric neighbourhood $U$ of the identity such that the $\sigma$-compact open subgroup $H_U:= \bigcup_{n \in \N} U^n$ is noncompact, then we are done, so assume that any such $H_U$ is compact. If we have chosen compact neighbourhoods $U_1, \ldots, U_k$ of the identity such that $H_{U_1} \lneq H_{U_2} \lneq \cdots \lneq H_{U_k}$, then, since $G$ is noncompact, we can choose a compact symmetric set $U_{k+1}$ that strictly contains $H_{U_k}$; so $H_{U_{k+1}} \gneq H_{U_k}$. We obtain a strictly increasing sequence of compact open subgroups of $G$, $(H_{U_n})_n$, and conclude that the $\sigma$-compact open subgroup $H:= \bigcup_{n \in\N}H_{U_n}$, which is covered by the open subgroups $H_{U_n}$, is noncompact.   
\end{proof}

The following proposition includes Theorem \ref{MainThm - metrizable case}. Recall that a locally compact group is compact-by-discrete if it has a compact open normal subgroup.

\bp  \rm \label{General Non-LI Prop} Let $G$ be a tdlc group satisfying one of the following conditions: \bi 
\item[(i)] $G$ is $\sigma$-compact and is not compact-by-discrete;

 \item[(ii)]  $G$ is noncompact and does not contain a  noncompact, $\sigma$-compact, compact-by-discrete  open subgroup;
 
 \item[(iii)] $G$ is noncompact and contains a sequence of neighbourhoods of the identity, $(U_n)_n$ , such that $\bigcap U_n$ has empty interior. 
 \ei  Then $\fB$ is not left introverted. 
 \ep

 \begin{proof}  If $G$ satisfies condition (ii) or (iii), take $H$ to be a noncompact, $\sigma$-compact open subgroup of $G$ (via Lemma \ref{Lemma Sigma-compact open subgroup}); if $G$ satisfies condition (i), let $H=G$. By Proposition \ref{Prop LI Quotients Open Subgroups}, it suffices to show that $\fB(H)$ is not left introverted.  Using the Kakutani-Kodaira Theorem \cite[Theorem 8.7]{HR}, take $N$ to be a compact normal subgroup of $H$ such that $H/N$ is metrizable; in the case that $G$ satisfies condition (iii), assume as well that  $N$ is contained in $\bigcap U_n$. In any case, $N$ is not open,  so $H/N$ is a first countable, $\sigma$-compact, noncompact, nondiscrete tdlc group, so  $\fB(H/N)$ is not left introverted by Theorem \ref{MainThm - metrizable case}.  By Proposition \ref{Prop LI Quotients Open Subgroups}, $\fB(H)$ is not left introverted, as needed.   
  \end{proof} 
 
 \br \rm 
 \noindent 1. If $G$ is a $\sigma$-compact locally compact group, then there is a sequence of neighbourhoods $(U_n)$ of the identity such that $\bigcap U_n = \{e_G\}$ (which has empty interior if $G$ is nondiscrete) if and only if $G$ is metrizable. 
 
 \smallskip 
 
 \noindent 2. Let  $G$ be a locally compact group with connected component of the identity, $G_e$.  If $G/G_e$ satisfies one of the three conditions stated in Proposition \ref{General Non-LI Prop}, then $\fB$ is not left introverted by Proposition \ref{Prop LI Quotients Open Subgroups}.  \er

 \br  \rm The theory of semigroup compactifications is commonly developed in the setting of semitopological semigroups \cite{Ber-Jun-Mil}. The main concepts from \cite{Ste-Sto} extend to this context as follows:  
  
 Let $S$ be a semitopological semigroup. For a subset $B$ of $S$ and $s \in S$, let  $$\text{$B\c s = \ell_s^{-1}(B)$ and $s \c B= r_s^{-1}(B)$, where $\ell_s, r_s: S \ra S: t \mapsto st, ts$.}$$    
The definitions of left (right) introversion of a $S$-invariant Boolean algebra $\B$ of clopen subsets of $S$, and left (right) Arens product on $\nu_\B S$, are now defined exactly as  before, and it is not difficult to check that Lemma 3.2 and Theorem 3.3 of \cite{Ste-Sto} remain true, following  nearly verbatim proofs; one only needs to note that  ``homomorphism" should replace ``isomorphism" in the statement of Lemma 3.2(i) and a few modest notational changes are required in the proof of Theorem 3.3(ii). Notice as well  that if $S$, a semitopological semigroup with at least one left identity $e$, replaces $G$ in  Proposition \ref{General Non-Introversion Prop}, then the proposition still  holds, with the same proof. Using this generalized version of Proposition \ref{General Non-Introversion Prop}, we conclude with another example.

Consider the Sorgenfrey line, $ \R_l$,  a zero dimensional, first countable, non-locally compact, non-$\sigma$-compact Hausdorff space. Observe that $(\R_l, +)$ is a topological semigroup (addition is jointly continuous on $\R_l$), but not a topological group (e.g., $\lim  1/n = 0$ but $\lim (-1/n) \neq 0$ in $\R_l$). To see that $\fB(\R_l)$ is not left introverted, let $B = \bigcup_{n\geq 2} [n, n+1/n)$, a clopen subset of $\R_l$, and let $a_n = n$, $c_k=1/k$. Then $(a_n)_{n\geq2}$ is contained in $B$, $\lim c_k = 0$, but $c_k + a_n \notin B$ for $n \geq k$, so condition (N$_B$) is satisfied.

 \er

\noindent {\sc Department of Mathematics and Statistics, University
of Winnipeg, 515 Portage Avenue, Winnipeg, MB, Canada, R3B 2E9}\\ {\sc Email: } {\tt  basmanmonterrubio-j@webmail.uwinnipeg.ca,  czyzowicz-t@webmail.uwinnipeg.ca, r.stokke@uwinnipeg.ca, thevenot-e@webmail.uwinnipeg.ca}

\end{document}